\documentclass[11pt]{amsart}

\usepackage{amsfonts,amsmath,latexsym,amssymb} 
\usepackage[T2A]{fontenc}
\usepackage[utf8]{inputenc}
\usepackage[english]{babel}
\usepackage{cite}

\newtheorem{theorem}{Theorem}
\newtheorem{lemma}{Lemma}
\newtheorem{remark}{Remark}

\newtheorem{corollary}{Corollary}

\newcommand{\RR}{\mathbb R}

\allowdisplaybreaks

\begin{document}
\title[Nagy type inequalities in metric measure spaces]{Nagy type inequalities in metric measure spaces and some applications
}

\author[V.~F.~Babenko]{Vladyslav Babenko}
\address{Department of Mathematical Analysis and Optimization, Oles Honchar Dnipro National University, Dnipro, Ukraine}

\author[V.~V.~Babenko]{Vira Babenko}
\address{Department of Mathematics and Computer Science, Drake University, Des Moines, USA}

\author[O.~V.~Kovalenko]{Oleg Kovalenko}
\address{Department of Mathematical Analysis and Optimization, Oles Honchar Dnipro National University, Dnipro, Ukraine}

\author[N.~V.~Parfinovych]{Nataliia Parfinovych}
\address{Department of Mathematical Analysis and Optimization, Oles Honchar Dnipro National University, Dnipro, Ukraine}
\begin{abstract}
We obtain a sharp Nagy type inequality in a metric space $(X,\rho)$ with measure $\mu$ that estimates the uniform norm of a function using its $\|\cdot\|_{H^\omega}$ -- norm determined by a modulus of continuity $\omega$, and a seminorm that is defined on a space of locally integrable  functions. 
We consider charges $\nu$ that are defined on the set of $\mu$--measurable subsets of $X$ and are absolutely continuous with respect to $\mu$. Using the obtained Nagy type inequality, we prove a sharp Landau -- Kolmogorov type inequality that estimates the uniform norm of a Radon -- Nikodym derivative of a charge via a $\|\cdot\|_{H^\omega}$ -- norm of this derivative, and a seminorm defined on the space of such charges. We also prove a sharp inequality for a hypersingular integral operator.
In the case $X=\RR_+^m\times \RR^{d-m}$, $0\le m\le d$, we obtain  inequalities  that estimate the uniform norm of a mixed derivative of a function using the uniform norm of the function and the $\|\cdot\|_{H^\omega}$ -- norm of its mixed derivative.
\end{abstract}

\keywords{
Nagy and Landau-Kolmogorov type inequalities, Stechkin's problem, charges, modulus of continuity, mixed derivative}


%% MSC codes here, in the form: 
\subjclass[2020]{ 26D10, 41A17, 41A44, 41A55 } 
\maketitle

\section{Introduction}
Inequalities that estimate norms of the intermediate derivatives of univariate or multivariate functions using the norms of the functions and their derivatives of higher order play an important role in many branches of Analysis and its applications. It appears that the richest applications are obtained from sharp inequalities of this kind, which attracts much interest to the inequalities with the smallest possible constants. For univariate functions, the results by Landau~\cite{Landau}, Kolmogorov~\cite{Kolmogorov1939}, and Nagy~\cite{Nagy}, are among the brightest ones in this topic.   A survey on the results for univariate and multivariate functions for the case of derivatives of integer and fractional order and further references can be found in~\cite{Arestov1996, BKKP, Babenko2012, Babenko22}. Article~\cite{Babenko05} (see also~\cite[Chapter 2]{BKKP}) contains a periodic analogue of Nagy's inequality;
% Connection {\color{blue} Тут не стоит писать Connection. В упомянутой ниже главе 2 приведены неравенство Надя и его аналог для периодических функций одного переменного. Так наверное и надо написать. И можно дать ссылку по этому поводу на мою с Пичуговым и Кофановым оригинальную статью, а потом на книгу. И соответствующий текст я переставляю немного ниже.}
% between the Landau -- Kolmogorov type  and  the Nagy type inequalities is discussed in~\cite[Chapter 2]{BKKP}
some recent results on the Nagy type inequalities are contained in~\cite{Kofanov}. Some inequalities of Landau -- Kolmogorov type for  Radon -- Nikodym derivative of charges defined on Lebesgue measurable subsets of an open cone  $C\subset \RR^d$ that are absolutely continuous with respect to the Lebesgue measure were obtained in~\cite{Babenko23}. In this article we study functional classes that are defined in terms of a majorant for modulus of continuity of the functions. The moduli of continuity  as independent functions and classes of functions with 
 given majorants of moduli of continuity of functions or of their derivatives  were introduced by Nikol'skii in~\cite{Nikolsky46}. Such classes were studied by many authors, see e.g.,~\cite[Chapter~7]{ExactConstants}.
 Extremal problems for such classes of non-real valued functions were considered in~\cite{Babenko16,VeraBabenko_JANO,Kovalenko20,Babenko21,Kriachko,BabenkoArxiv,Babenko22b}. Extremal problems for various hyper-singular integral operators on classes of univariate and multivariate functions defined by a majorant on their modulus of continuity were considered in~\cite{Babenko07,babenko2010,BPP2010,Parfinovych}.

Let $X$, $Y$, and $Z$ be linear spaces equipped with seminorm $\|\cdot\|_X$, norm $\|\cdot\|_Y$, and seminorm $\|\cdot\|_Z$ respectively. A linear operator $A\colon X\to Y$ is called bounded, if 
\[
\| A\|=\| A\|_{X\to Y}=\sup_{\| x\|_X\le 1}\| Ax\|_Y<\infty.
\]
Otherwise the operator $A$ is called unbounded.
By $\mathcal{L}(X,Y)$ we denote  the space of all linear bounded operators $S\colon X\to Y$. 

Let $A\colon X \to Y\cap Z$ be a homogeneous operator (not necessarily linear) with the domain $D_A \subset X$. Let also $\mathfrak{M}=\{ x\in D(A) \colon \| Ax\|_Z\le1\}$.
For the operator $A$ and an operator  $S\in \mathcal{L}(X,Y)$ we set
$$
    U(A,S; \mathfrak{M})=\sup\{\|Ax-Sx\|_Y\colon x\in \mathfrak{M}\}.
$$
Note that for each $x\in D_A$ one has
\[
\| Ax-Sx\|_Y\le U(A,S; \mathfrak{M})\| Ax\|_Z.
\]

The Stechkin problem of approximation of a generally speaking unbounded operator by linear bounded operators on the class $\mathfrak{M}$ is formulated as follows. For a given number $N$ find the quantity
\begin{equation}\label{bestappr}
    E_N(A,\mathfrak{M})=\inf\left\{ U(A,S;\mathfrak{M})\colon S\in{\mathcal L}(X,Y), \|S\|\le N\right\}.
\end{equation}  The statement of a somewhat more general problem, first important results, and
solutions to this problem for differential operators of small orders were presented in~\cite{Stechkin1967}. For a survey of further results on this problem see~\cite{ Arestov1996}. 
The Stechkin problem, in turn, is intimately connected to Landau-Kolmogorov type inequalities. The following well-known theorem (which we formulate in a convenient for us form) describes this connection. 

\begin{theorem} For any $x\in D_A$ and arbitrary $S\in \mathcal{L}(X,Y)$ the following Landau -- Kolmogorov -- Nagy type inequality holds
\begin{equation}\label{abstraddineq} 
 \| Ax\|_Y\le  \| Ax-Sx\|_Y+\| S\|\|x\|_X
\le U(A,S;\mathfrak{M})\| Ax\|_Z+\|S\|\cdot \|x\|_X,    
\end{equation}
    and, therefore,
$$
     \forall x\in D_A,\,\forall N>0, \, \|Ax\|_Y\le E_N(A,\mathfrak{M})\| Ax\|_Z+N\|x\|_X.  
$$
If in addition there exist $\overline{S}\in \mathcal{L}(X,Y)$ and  $\overline{x}\in \mathfrak{M}$ such that
$$
    \|A\overline{x}\|_Y=\| A\overline{x}-\overline{S}\overline{x}\|_Y+\| \overline{S}\|\|\overline{x}\|_X=U(A,\overline{S};\mathfrak{M})+\|\overline{S}\|\cdot \|\overline{x}\|_X,
$$
then 
$$
E_{\left\|\overline{S}\right\|}(A, \mathfrak{M})=U(A,\overline{S};\mathfrak{M})=\| A\overline{x}\|-\left\|\overline{S}\right\|\|\overline{x}\|_X,
$$
and the operator $\overline{S}$ is optimal for problem~\eqref{bestappr}.
\end{theorem}

\begin{remark}
In Stechkin's article~\cite{Stechkin1967} it is assumed that  $X$ and $Y$ are Banach spaces. However, as it is easy to see, completeness and even presence of a norm in $X$ is not necessary. It is sufficient to have a seminorm in $X$. Completeness of $Y$ is also not necessary.
\end{remark}

In Section~\ref{sec::notations} we give necessary notations and definitions. In Section~\ref{s::NagyTypeInequality} we obtain a sharp Nagy type inequality in a metric space $(X,\rho)$ with measure $\mu$ that estimates the uniform norm of a function using its $\|\cdot\|_{H^\omega}$ -- norm determined by a modulus of continuity $\omega$, and a seminorm that is defined on a space of locally integrable  functions. 
 In Section~\ref{s::metricSobolevSpaces} we prove a sharp inequality of Nagy type in the context of metric Sobolev spaces.
Using the inequality from Section~\ref{s::NagyTypeInequality}, in Section~\ref{s::LKforCharges} we prove a sharp Landau -- Kolmogorov type inequality that estimates the uniform norm of a Radon -- Nikodym derivative of a charge from a particular class of charges via a $\|\cdot\|_{H^\omega}$ -- norm of this derivative, and a seminorm defined on the space of charges.
In  Section~\ref{s::hypersingularIntegral} we obtain a sharp Landau -- Kolmogorov type inequality for generalized hypersingular operators.
Finally, in Section~\ref{s::mixedDerivative} we suppose that $X=\RR_+^m\times \RR^{d-m}$, $0\le m\le d$, and obtain  inequalities  that estimate the uniform norm of a mixed derivative of a function $f\colon X\to \RR$ using the uniform norm of the function and the norm of its mixed derivative which is defined with the help of some modulus of continuity. 

We use the following scheme to obtain the main results of the article. We define an appropriate bounded operator $S$, give an estimate for the quantity $U(A,S,\mathfrak{M})$, and plugging it into~\eqref{abstraddineq}, we obtain a Nagy or Landau -- Kolmogorov type inequality. Then we prove its sharpness. Theorem~1 shows that we simultaneously obtain a solution to the corresponding Stechkin problem.

\section{Notations and definitions}\label{sec::notations}

Let $(X,\rho)$ be a metric space with a Borel measure $\mu$. Assume that $X$ is a  commutative monoid (i.e., an associative and commutative  binary operation $+$ is defined on $X$, and there exists an element $\theta\in X$ such that $x  + \theta = \theta + x = x$ for all $x\in X$) such that for each measurable set $Q\subset X$ and each $x\in X$ one has
\[
\mu (x+Q)=\mu (Q).
\]
Suppose that for all $x,y\in X$, 
\[
\rho(x+y,x)\le \rho(y,\theta).
\]
Everywhere below $B_h=B_h(\theta)$ is an open ball of radius $h>0$ with center $\theta$; we suppose that $0<\mu(B_h)<\infty$ and $B_h\neq \{\theta\}$ for all $h>0$.

An invariant Haar measure on a locally compact group with metrical topology, see e.g.~\cite{Loomis}, is an important example.

For a measurable set $Q\subset X$ by $L_1(Q)$ ($L_\infty(Q)$) we denote the space of functions $f\colon Q\to\RR$ integrable (resp. essentially bounded) on $Q$  with the corresponding norm. 
By $L_{{\rm loc}}(X)$ we denote the space of all functions $f\colon X\to\RR$ that are integrable on each open ball of $X$. In the space  $L_{\rm loc}(X)$ we introduce a family of seminorms
\[
\rfloor f\lceil_h
=
\sup_{x\in X}\left|\;\int_{x+B_h}f(u)d\mu(u)\right|,\, h>0,
\]
and a seminorm
\[
\rfloor f\lceil=\sup_{h>0}\rfloor f\lceil_h.
\]
By $L_{\rfloor \cdot\lceil_h}(X)$ ($L_{\rfloor \cdot\lceil}(X)$) we denote the family of functions $f\in L_{\rm loc}(X)$ with a finite seminorm $\rfloor \cdot\lceil_h$ (resp. $\rfloor \cdot\lceil$). It is clear that the space  $L_1(X)$ is contained in each of these sets.

By $C(X)$ we denote the space of all continuous functions $f\colon X\to\RR$; by  $C_b(X)$ the space of functions  $f\in C(X)$ with a finite norm 
\[
\| f \|_{C(X)}=\sup_{x\in X}|f(x)|;
\]
by $\mathcal{B}(X)$ we denote the space of bounded functions $f\colon X\to \RR$ with a norm
\[
\| f \|_{\mathcal{B}(X)}=\sup_{x\in X}|f(x)|.
\]
%{\color{red}We also use the standard notations $L_p(X)$, $p\in [1,\infty]$.}
Everywhere below we assume that the measure  $\mu$ is such that $C(X)\subset L_{\rm loc}(X)$.

% {\color{red} Ниже используются также обозначения $L_\infty(X)$ and $L_1(X)$. Стоит наверное сказать об этом. Дескать наряду с введенными мы используем также эти стандартные обозначения}

Let $\omega$ be a modulus of continuity i.e., a non-negative, non-decreasing, semi-additive function $\omega\colon [0,\infty)\to[0,\infty)$ such that $\omega(0)=0$.  By $H^\omega(X)$ we denote the space of functions $f\colon X\to \RR$ such that 
\[
\| f\|_{H^\omega(X)} :=\sup_{x,y\in X, x\neq y}\frac{|f(x)-f(y)|}{\omega(\rho(x.y))}<\infty.
\]

\section{A Nagy type inequality}\label{s::NagyTypeInequality}
For each $h>0$ we define an operator 
$S_h\colon L_{\rfloor \cdot\lceil_h}(X)\to \mathcal{B}(X)$ by the following rule:
\[
S_hf(x)=\frac 1{\mu(B_h)}\int_{B_h}f(x+u)d\mu(u).
\]
It is clear that this operator is bounded and 
\begin{equation}\label{normSh}
 \| S_h\|_{L_{\rfloor \cdot\lceil_h}(X)\to \mathcal{B}(X)}=\frac 1{\mu(B_h)}.   
\end{equation}

We need the following result, which  is sometimes called an Ostrowski type inequality.  A rather general version of such kind of results is contained in~\cite[Theorem~2]{Babenko22b}.
\begin{lemma}\label{l::ostrowski}
 If $f\in H^\omega(X)$, then for each $h>0$ one has 
\begin{equation}\label{Ostrowskii}
 \|f(\cdot)-S_hf(\cdot)\|_{\mathcal{B}(X)}\le\frac {\| f\|_{H^\omega(X)}}{\mu(B_h)}\int_{B_h}\omega(\rho (u,\theta))d\mu(u).
 \end{equation}
 Inequality~\eqref{Ostrowskii} is sharp and becomes equality for all functions 
 \[
 f_\omega(x)=c\pm \omega(\rho (u,\theta)), c\in \RR.
 \]
\end{lemma}
\begin{proof}
For each $x\in X$ we have
\begin{multline*}
|f(x)-S_hf(x)|=\left|f(x)-\frac 1{\mu(B_h)}\int_{B_h}f(x+u)d\mu(u)\right|
\\ \le
\frac 1{\mu(B_h)}\int_{B_h}|f(x)-f(x+u)|d\mu(u)
 \le
\frac{\| f\|_{H^\omega(X)}}{\mu(B_h)}\int_{B_h}\omega(\rho (x+u,x))d\mu(u)
\\ \le
\frac{\| f\|_{H^\omega(X)}}{\mu(B_h)}\int_{B_h}\omega(\rho (u,\theta))d\mu(u)
\end{multline*}
and inequality~\eqref{Ostrowskii} is proved. For the function  $f_\omega$ one has 
\begin{equation}\label{fOmegaNorm}
\| f_\omega\|_{H^\omega(X)}=1.
\end{equation}
Indeed, for all $x,y\in X$,
$$
|f_\omega(x)- f_\omega(y)| = |\omega(\rho(x,\theta)) - \omega(\rho(y,\theta))|\leq \omega(|\rho(x,\theta) -\rho(y,\theta)|)\leq \omega(\rho(x,y)),
$$
hence $\|f_\omega\|_{H^\omega(X)}\leq 1$; since there exists $y\in B_h\setminus\{\theta\}$, we obtain 
$$
\sup_{x\in X, x\neq \theta} |f_\omega(x) - f_\omega(\theta)|\geq |f_\omega(y) - f_\omega(\theta)| = \omega(\rho(y,\theta)),
$$
which implies $\|f_\omega\|_{H^\omega(X)}\geq 1$, and hence~\eqref{fOmegaNorm}.
Moreover,
\[
\|f(\cdot)-S_hf(\cdot)\|_{\mathcal{B}(X)}
\ge 
|f(\theta)-S_hf(\theta)| 
=
\frac {1}{\mu(B_h)}\int_{B_h}\omega(\rho (u,\theta))d\mu(u),
\]
which together with~\eqref{fOmegaNorm} implies that inequality~\eqref{Ostrowskii} is sharp.
\end{proof}

The following theorem contains a variant of the Nagy type inequality. For $\alpha\in \RR$ we set $\alpha_+ :=\max\{\alpha,0\}$.

\begin{theorem}\label{th::HwL-K}
    If $h>0$ and $f\in H^\omega(X)\cap  L_{\rfloor \cdot\lceil_h}(X)$, then
\begin{multline}\label{Nagy_h}
\|f\|_{\mathcal{B}(X)}\le \| f-S_hf\|_{\mathcal{B}(X)}+\| S_h\|_{L_{\rfloor \cdot\lceil_h}(X)\to \mathcal{B}(X)}\| f\|_{L_{\rfloor \cdot\lceil_h}(X)}
\\ \le
\frac {\| f\|_{H^\omega(X)}}{\mu(B_h)}\int_{B_h}\omega(\rho (u,\theta))d\mu(u)+\frac {\rfloor f\lceil_h}{\mu(B_h)}.   
\end{multline}
The inequality is sharp and turns into equality for the function
\begin{equation}\label{f_eh}
f_{e,h}(x)=(\omega(h)-\omega(\rho(x,\theta)))_+.
\end{equation}
Moreover, $f_{e,h}\in H^\omega(X)\cap  L_{\rfloor \cdot\lceil}(X)$, $\rfloor f_{e,h}\lceil = \rfloor f_{e,h}\lceil_h$, and hence for each $h>0$ the inequality 
\begin{equation}\label{Nagy}
\|f\|_{\mathcal{B}(X)}
\le 
\frac {\| f\|_{H^\omega(X)}}{\mu(B_h)}\int_{B_h}\omega(\rho (u,\theta))d\mu(u)+\frac {\rfloor f\lceil }{\mu(B_h)}  
\end{equation}
holds and is sharp on the class $ H^\omega(X)\cap  L_{\rfloor \cdot\lceil}(X)$.
\end{theorem}

\begin{remark}
In the case, when $\omega(t)=t^\alpha$, $0<\alpha\le 1$, and the space  $(X,\rho,\mu)$ satisfies the following $s$-regularity type condition   
    $$
    \exists b>0, \exists s>0 \colon \forall h>0,\mu(B_h)\ge bh^s,
    $$ 
inequality~\eqref{Nagy} can be written in a multiplicative form. We do not adduce the details.
\end{remark}

\begin{proof}
For each  $x\in X$, due to Lemma~\ref{l::ostrowski}, the definition of the operator  $S_h$, and equality~\eqref{normSh},
\begin{multline*}
|f(x)|\le |f(x)-S_hf(x)|+|S_hf(x)|
\\\le 
\|f(x)-S_hf(x)\|_{B(X)}+\|S_h\|_{L_{\rfloor \cdot\lceil_h}(X)\to B(X)}\rfloor f\lceil_h
\\ \le 
\frac {\| f\|_{H^\omega(X)}}{\mu(B_h)}\int_{B_h}\omega(\rho (u,\theta))d\mu(u)+\frac {\rfloor f\lceil_h}{\mu(B_h)},
\end{multline*}
which implies inequality~\eqref{Nagy_h}. Inequality~\eqref{Nagy} is a consequence of inequality~\eqref{Nagy_h}.

For the function $f_{e,h}$ we have 
$\| f_{e,h}\|_{\mathcal{B}(X)}= f_{e,h}(\theta) = \omega(h)$, $\| f_{e,h}\|_{H^\omega(X)}=1$, and 
$$
\rfloor f_{e,h}\lceil_h=\omega(h)\mu(B_h)-\int_{B_h}\omega(\rho (u,\theta))d\mu(u).
$$
Indeed, on the one hand, 
$$\rfloor f_{e,h}\lceil_h\geq \int_{\theta+B_h}f_{e,h}(u)d\mu(u) = \omega(h)\mu(B_h)-\int_{B_h}\omega(\rho (u,\theta))d\mu(u),$$
and on the other hand, for each $x\in X$, due to monotonicity of $\omega$,
\begin{multline*}
\int_{x+B_h}f_{e,h}(u)d\mu(u) 
=
\int_{(x+B_h)\cap B_h}f_{e,h}(u)d\mu(u)
\\ \leq 
\int_{B_h}f_{e,h}(u)d\mu(u) = \omega(h)\mu(B_h)-\int_{B_h}\omega(\rho (u,\theta))d\mu(u).
\end{multline*}
Direct computations now show that inequality~\eqref{Nagy_h} becomes equality on the function $f_{e,h}$.

Finally, the same arguments as during computation of the quantity $\rfloor f_{e,h}\lceil_h$ show that $\rfloor f_{e,h}\lceil = \rfloor f_{e,h}\lceil_h$, and hence inequality~\eqref{Nagy} is also sharp.
\end{proof}

\begin{corollary}\label{th::HwCorollary}
        If $h>0$ and $f\in H^\omega(X)\cap  L_{1}(X)$, then
$$
\|f(x)\|_{\mathcal{B}(X)}\le \frac {\| f\|_{H^\omega(X)}}{\mu(B_h)}\int_{B_h}\omega(\rho (u,\theta))d\mu(u)+\frac 1{\mu(B_h)}\left\|f\right\|_{L_1(X)}.
$$
The inequality is sharp. It becomes equality on the function $f_{e,h}$ defined by~\eqref{f_eh}.
\end{corollary}
\begin{proof}
    Since $L_1(X)\subset L_{\rfloor \cdot\lceil_h}(X)$ and $\|f_{e,h}\|_{L_1(X)} = \rfloor f_{e,h}\lceil_h$, the statement of the corollary follows from Theorem~\ref{th::HwL-K}.
\end{proof}

\section{Nagy type inequalities in metric Sobolev spaces}\label{s::metricSobolevSpaces}
In this section we consider metric Sobolev spaces with essentially bounded upper gradients, which are defined as follows.
% {\color{red} Я бы начал примерно такой фразой. Здесь мы будем рассматривать метрические Соболевские пространства с существенно ограниченными верхними градиентами, которые определяются следующим образом.}
Let $(X,\rho_X)$ and $(Y,\rho_Y)$ be metric spaces such as in the previous sections. For a modulus of continuity $\omega$ we define the space $W^{1,\omega}(X,Y)$ as the space of all functions $f\colon X\to Y$ with the following property (cf.~\cite[Chapter~5]{ambrosio2004topics} and~\cite[Chapter~10.2]{heinonen2015sobolev}). There exists a non-negative function $G=G_f\in L_\infty(X)$ and a set $N = N_f\subset X$  such that $\mu (N)=0$ and
\begin{equation}\label{upperGradient}  
\rho_Y(f(x),f(y))\le (G(x)+G(y))\cdot \omega(\rho_X(x,y))\text{ for all } x,y\in X\setminus N.
\end{equation}
We call $G$ an upper gradient of $f$.
%{\color{red} Здесь,а может и в самом начале можно написать следующее. Конечно, наряду с функцией $\rho_X(x,y)$ функция $\omega(\rho_X(x,y))$ также является (деформированной) метрикой, так что определенное выше Пространство можно рассматривать как деформацию стандартного метрического Соболевского пространства. По видимому одной из первых работ, посвященных Дефомациям метрических просттранств являетcа работа~\cite{Timan}. Статью я посылаю на почту. }

Let $(Y,\rho_Y)$ be the space of reals with the usual metric. The technique developed in the proof of the previous theorem allows to prove the following result, which thus can be  in some sense considered as a corollary of Theorem~\ref{th::HwL-K}.
% {\color{red} Здесь возможно стоит сказать, что техника, развитая при доказательстве предыдущей теоремы, позволяет доказать теорему 3, иак что ее можно рассматривать в некотором роде как приложение предыдущей теоремы. И в глае про заряды тоже можно привести без доказательства следствие для метрических Соблолевских пространств}
\begin{theorem}\label{th::MSS}
    Let $f\in W^{1,\omega}(X,\RR)\cap  L_{\rfloor \cdot\lceil_h}(X)$ and $G_f$ be an upper gradient of $f$. Then for any $h>0$ the following inequality holds
    \begin{equation}\label{MSS}
        \| f\|_{L_\infty(X)}\le \frac{2\| G_f\|_{L_\infty(X)}}{\mu(B_h)}\int_{B_h}\omega(\rho_X(\theta,u) )d\mu(u)+\frac {\rfloor f\lceil_h}{\mu(B_h)}.
    \end{equation}
    Inequality~\eqref{MSS} is sharp in the sense that there exists a function $f$ and its upper gradient $G_{f}$ for which the inequality becomes equality.
\end{theorem}
\begin{proof}
For almost all $x\in X$, 
\begin{gather*}
    |f(x)|
    \leq 
    \left|f(x) - \frac{1}{\mu(B_h)}\int_{x+B_h}f(y)d\mu(y)\right| +  \frac{1}{\mu(B_h)}\left|\int_{x+B_h}f(y)d\mu(y)\right|
    \\ \leq
    \frac{1}{\mu(B_h)}\int_{x+B_h}\left|f(x) - f(y)\right|d\mu(y) +  \frac {\rfloor f\lceil_h}{\mu(B_h)}
    \\ \leq 
    \frac{1}{\mu(B_h)}\int_{B_h}\left|G_f(\theta) - G_f(u)\right|\omega(\rho_X(\theta, u))d\mu(u) +  \frac {\rfloor f\lceil_h}{\mu(B_h)}
    \\ \leq 
        \frac{2\| G_f\|_{L_\infty(X)}}{\mu(B_h)}\int_{B_h}\omega(\rho_X(\theta,u) )d\mu(u)+  \frac {\rfloor f\lceil_h}{\mu(B_h)},
\end{gather*}
and inequality~\eqref{MSS} is proved. The inequality becomes equality for the function $f = f_{e,h}$ defined in~\eqref{f_eh} and its upper gradient $G \equiv \frac 12$. Indeed, the fact that inequality~\eqref{MSS} becomes equality for such functions $f$ and $G$ can be verified directly; the fact that $G$ is indeed an upper gradient for $f$ (with $N = \emptyset$) can be proved as follows. Let $x,y\in X$, $u = \rho_X(x,\theta)$ and $v  =\rho_X(y,\theta)$. We can assume that $u\leq v$. If $h\leq u\leq v$, then $f(x) = f(y) = 0$ and inequality~\eqref{upperGradient} holds. If $0\leq u<  h \leq v$, then
\begin{multline*}
    |f(x) - f(y)| 
    = 
    \omega(h) - \omega (u)
    \leq 
    \omega(h-u)\leq \omega (v-u) 
    \\= 
    \omega(\rho_X(y,\theta) - \rho_X(x,\theta))
     \leq 
    \omega(\rho_X(x,y)).
\end{multline*}
The case $0\leq u\leq v\leq h$ can be considered similarly.
\end{proof}
%{\color{red}
Similarly to Corollary~\ref{th::HwCorollary} one can prove the following result.
\begin{corollary}
    Let $f\in W^{1,\omega}(X,\RR) \cap  L_{1}(X)$ and $G_f$ be an upper gradient of $f$. Then for any $h>0$ the following inequality holds
    $$
        \| f\|_{L_\infty(X)}\le \frac{2\| G_f\|_{L_\infty(X)}}{\mu(B_h)}\int_{B_h}\omega(\rho_X(\theta,u) )d\mu(u)+\frac {\| f\|_{L_1(X)}}{\mu(B_h)}.
    $$
    The inequality is sharp in the sense that there exists a function $f$ and its upper gradient $G_{f}$ for which the inequality becomes equality.
\end{corollary}
%}
\section{Landau -- Kolmogorov type inequalities for charges}\label{s::LKforCharges}
By $\mathfrak{N}(X)$ we denote the family of charges $\nu$ defined on the family of all $\mu$--measurable subsets of $X$ and that are absolutely continuous with respect to the measure $\mu$, see e.g.,~\cite[Chapter 5]{Berezanski}. By the Radon--Nikodym theorem, for a charge $\nu\in \mathfrak{N}(X)$ there exists an integrable function $f\colon X\to\RR$ such that for arbitrary measurable set $Q\subset X$
\begin{equation}\label{ac_charge}
\nu(Q)=\int_Qf(x)d\mu(x).
\end{equation}
This function $f$ is called the Radon--Nikodym derivative of the charge $\nu$ with respect to the measure $\mu$ and will be denoted by $D_\mu\nu$.
The family $\mathfrak{N}(X)$ is a linear space with respect to the standard addition and multiplication by a real number. Define a family of seminorms $\{\rceil \cdot\lfloor_h,\; h>0\}$ as follows:
\[
\rceil \nu\lfloor_{h}=\|\nu(\cdot +B_h)\|_{B(X)}.
\]
It is clear that if a charge $\nu$ and a function $f$ are related via~\eqref{ac_charge}, then
\[
\rceil \nu\lfloor_h=\rfloor f\lceil_h.
\]
For $h> 0$ by $\mathfrak{N}_{\rceil \cdot\lfloor_h}(X)$ we denote the set of charges  $\nu\in\mathfrak{N}(X)$ with a finite seminorm $\rceil \cdot\lfloor_h$. 
%The set of charges $\nu\in\mathfrak{N}(X)$ such that $D_\mu\nu\in H^\omega(X)$ will be denoted by  $\mathfrak{N}^{D_\mu,\omega}(X)$. 
\begin{theorem}
If $h>0$ and $\nu\in \mathfrak{N}_{\rceil \cdot\lfloor_h}(X)$ is such that $D_\mu\nu\in H^\omega(X)$, then
\begin{multline*}
   \left\| D_\mu\nu\right\|_{\mathcal{B}(X)}
   \le
   \left\| D_\mu\nu-\overline{S}_h\nu\right\|_{\mathcal{B}(X)}+\| \overline{S}_h\|_{\mathfrak{N}_{\rceil \cdot\lfloor_h}(X)\to \mathcal{B}(X)}\rceil \nu\lfloor_h
    \\\le 
   \frac {\|  D_\mu\nu\|_{H^\omega(X)}}{\mu(B_h)}\int_{B_h}\omega(\rho (u,\theta))d\mu(u)+\frac {\rceil \nu\lfloor_h}{\mu(B_h)},
\end{multline*}
where the operator $\overline{S}_h\colon \mathfrak{N}_{\rceil \cdot\lfloor_h}(X)\to \mathcal{B}(X)$ is defined by
$$
\overline{S}_h\nu(x)=\frac{\nu(x+B_h)}{\mu(B_h)}.
$$
The inequality is sharp, and becomes an equality for the charge $\nu_{e,h}$ such that $D_\mu\nu_{e,h}=f_{e,h}$, where the function $f_{e,h}$ is defined by~\eqref{f_eh}.
\end{theorem}
\begin{remark}
This theorem generalizes Theorem~3 from~\cite{Babenko23}.    
\end{remark}
\begin{proof}
    It is enough to apply Theorem~\ref{th::HwL-K} to the function $f=D_\mu\nu$ and notice that after a change of variables in the integral, we obtain $\overline{S}_h\nu=S_h f.$
\end{proof}

\begin{corollary}
If $h>0$, $\nu\in \mathfrak{N}_{\rceil \cdot\lfloor_h}(X)$ is such that $D_\mu\nu\in W^{1,\omega}(X,\RR)$, and $G_\nu$ is an arbitrary upper gradient of $D_\mu \nu$, then the following inequality holds
    $$
        \| D_\mu\nu\|_{L_\infty(X)}\le \frac{2\| G_\nu\|_{L_\infty(X)}}{\mu(B_h)}\int_{B_h}\omega(\rho_X(\theta,u) )d\mu(u)+\frac {\rfloor \nu\lceil_h}{\mu(B_h)}.
   $$
   The inequality is sharp in the sense that there exists a charge $\nu$ and an upper gradient $G_{\nu}$ for which the inequality becomes equality.
\end{corollary}
\begin{proof}
    The theorem follows from Theorem~\ref{th::MSS}.
\end{proof}
\section{Inequalities for generalized hypersingular integrals}\label{s::hypersingularIntegral}
Let $P\colon \RR_+\to\RR_+$ be a locally integrable function such that for some $h>0$
\[
\int_{X\setminus B_h}P(\rho(u,\theta))d\mu(u)<\infty.
\]
Define the following operator, which can be considered as a hypersingular integral 
\[
\mathfrak{D}_Pf(x)=\int_X(f(x)-f(x+u))P(\rho(u,\theta))d\mu(u).
\]
We also consider the following truncated hypersingular integral
\[
\overline{\mathfrak{D}}_{P,h}f(x)=\int_{X\setminus B_h}(f(x)-f(x+u))P(\rho(u,\theta))d\mu(u).
\]
It is easy to see that 
\[
\|\overline{\mathfrak{D}}_{P,h}\|_{C_b(X)\to C_b(X)}={ 2}\int_{X\setminus B_h}P(\rho(u,\theta))d\mu(u).
\]
%{\color{red} Нам же еще нужно, чтоб $$
%\int_{X\setminus B_h}P(\rho(u,\theta))d\mu(u)<\infty,
%$$
%правильно? Стоит ли добавить такое требование в формулировку теореме? Мне кажется, что да - т.е. подчеркнуть требование, чтоб оба условия виполнялись для одного и того же $h$.}
\begin{theorem}
    Let $f\in H^\omega(X)\cap C_b(X)$. If for some $h>0$
$$
\int_{B_h}\omega(\rho(u,\theta))P(\rho(u,\theta))d\mu(u)<\infty\text{ and } \int_{X\setminus B_h}P(\rho(u,\theta))d\mu(u)<\infty,
$$
then
\begin{multline}\label{hiper}
\|\mathfrak{D}_Pf\|_{\mathcal{B}(X)}\le \|\mathfrak{D}_Pf-\overline{\mathfrak{D}}_{P,h}f\|_{\mathcal{B}(X)}+\|\overline{\mathfrak{D}}_{P,h}\|_{C_b(X)\to C_b(X)}\| f\|_{C_b(X)}
\\
\le  \| f\|_{H^\omega(X)}\int_{B_h}\omega(\rho(u,\theta))P(\rho(u,\theta))d\mu(u) \\+ 2\| f\|_{ C_b(X)}\int_{X\setminus B_h}P(\rho(u,\theta))d\mu(u).
\end{multline}
The inequality is sharp and turns into equality for the function 
\[
f_{e,\omega}(u)=\begin{cases}
    \omega(\rho(u,\theta))-\frac 12\omega(h), & \rho(u,\theta)\le h,\\ 
\frac 12\omega(h), & \rho(u,\theta)\ge h.
\end{cases}
\]
\end{theorem}
\begin{proof}
We have
\begin{gather*}
    \|\mathfrak{D}_Pf-\overline{\mathfrak{D}}_{P,h}f\|_{\mathcal{B}(X)} 
    =
    \sup_{x\in X}\left|\int_{B_h}(f(x)-f(x+u))P(\rho(u,\theta))d\mu(u)\right|
    \\ \leq 
    \sup_{x\in X}\| f\|_{H^\omega(X)}\int_{B_h}\omega(\rho(x,x+u))P(\rho(u,\theta))d\mu(u)
    \\ \leq 
    \|f\|_{H^\omega(X)}\int_{B_h}\omega(\rho(u,\theta))P(\rho(u,\theta))d\mu(u).
\end{gather*}
Thus
\begin{multline*}
\|\mathfrak{D}_Pf(x)\|_{\mathcal{B}(X)}\le \|\mathfrak{D}_Pf-\overline{\mathfrak{D}}_{P,h}f\|_{\mathcal{B}(X)}+\|\overline{\mathfrak{D}}_{P,h}\|_{C_b(X)\to C_b(X)}\| f\|_{C_b(X)}
\\ \le  
\|f\|_{H^\omega(X)}\int_{B_h}\omega(\rho(u,\theta))P(\rho(u,\theta))d\mu(u)
\\+
2\| f\|_{C_b(X)}\int_{X\setminus B_h}P(\rho(u,\theta))d\mu(u),
\end{multline*}
and inequality~\eqref{hiper} is proved.
For the function $f_{e,\omega}$ we have $\|f_{e,\omega}\|_{C_b(X)} = \frac 12\omega(h)$, $\|f_{e,\omega}\|_{H^\omega(X)} = 1$, and
\begin{multline*}
\|\mathfrak{D}_Pf_{e,\omega}\|_{\mathfrak{B}(X)} 
= 
-\mathfrak{D}_Pf_{e,\omega}(\theta) 
\\= \int_{B_h}\omega(\rho(u,\theta))P(\rho(u,\theta))d\mu(u)
+\omega(h) \int_{X\setminus B_h}P(\rho(u,\theta))d\mu(u),
\end{multline*}
thus the inequality becomes equality for the function $f_{e,\omega}$.
\end{proof}
\section{Inequalities for mixed derivatives}\label{s::mixedDerivative}
Assume that $X=\RR^d_{m,+}:=\RR^m_+\times \RR^{d-m}, 1\le m\le d$,  $\mu$ is the Lebesgue measure  in $\RR^d_{m,+}$, $\rho(x,y)=\max_{i=1,\ldots,d}|x_i-y_i|$, so that $B_h=(0,h)^m\times (-h,h)^{d-m}$. In this section for brevity we write $dx$ instead of $d\mu(x)$.
For a locally integrable function $f\colon X\to\RR$ set ${\bf I}=(1,\ldots,1)\in \RR^d$ and
$$\partial_{\bf I}f=\frac{\partial ^d f}{\partial x_1\ldots\partial x_d},$$
where the derivatives are understood in the distributional sense.

Let $\{e_i\}$ be the standard basis in $\RR^d$. For $i = 1,\ldots, d$ and $h>0$ we set
\[
\Delta^+_{i,h}f(x):=f(x+he_i)-f(x),
\text{ and }
\Delta_{i,h}f(x):=f(x+he_i)-f(x-he_i).
\]
In virtue of the Fubini theorem, for almost all $x\in \RR^d_{m,+}$ one has
\begin{equation}\label{mixdiff}
\int_{x+B_h} \partial_{\bf I}f(u)du=(\Delta^+_{1,h}\circ\ldots\circ\Delta^+_{m,h}\circ\Delta_{m+1,h}\circ\ldots\circ\Delta_{d,h})f(x).   \end{equation}
Define an operator $\mathfrak{S}_h\colon L_\infty(C)\to L_\infty(C)$, setting
\[
\mathfrak{S}_hf(x)=\frac 1{2^{d-m}h^d}\left(\Delta^+_{1,h}\circ\ldots\circ\Delta^+_{m,h}\circ\Delta_{m+1,h}\circ\ldots\circ\Delta_{d,h}\right)f(x), h > 0.
\]

\begin{theorem}
    If $h>0$ and $f\in \mathcal{B}(\RR^d_{m,+})$ is such that $\partial_{\bf{I}}f\in H^\omega(\RR^d_{m,+})$, then
\begin{multline}\label{partial_omega}
\|\partial_
{\bf{I}}f\|_{\mathcal{B}(\RR^d_{m,+})}
\le
\left\| \partial_{\bf I}f-\mathfrak{S}_hf\right\|_{\mathcal{B}(\RR^d_{m,+})}+\| \mathfrak{S}_h\|\|f\|_{\mathcal{B}(\RR^d_{m,+})}
\\\le 
\frac {\| \partial_
{\bf{I}}f\|_{H^\omega(\RR^d_{m,+})}}{2^{d-m}h^d}\int_{B_h}\omega(\rho (u,\theta))du +\frac {2^m}{h^d}\| f\|_{\mathcal{B}(\RR^d_{m,+})}.
\end{multline}

In the case, when $\omega(t) = t^\alpha$, $\alpha\in (0,1]$, the following multiplicative inequality holds:
\begin{equation}\label{partialOmegaMultiplicative}
\|\partial_
{\bf{I}}f\|_{\mathcal{B}(\RR^d_{m,+})}
\leq 
2^{\frac{m\alpha}{d+\alpha}}\left(\frac {d+\alpha}{\alpha}\right)^{\frac {\alpha}{d+\alpha}}\cdot \| f\|_{\mathcal{B}(\RR^d_{m,+})}^{\frac{\alpha}{d+\alpha}} \cdot\| \partial_
{\bf{I}}f\|_{H^\omega(\RR^d_{m,+})}^\frac{d}{d+\alpha} .
\end{equation}
For $m=0$ and $m = 1$ these inequalities are sharp.
\end{theorem}

\begin{remark}
This theorem generalizes Theorem~5 from~\cite{Babenko23}.    
\end{remark}

\begin{proof} 
Applying inequality~\eqref{Nagy_h} to $\partial_{\bf{I}}f$ and taking into account that $\mu(B_h)=2^{d-m}h^d$, we obtain
\begin{multline*}
\|\partial_
{\bf{I}}f\|_{\mathcal{B}(\RR^d_{m,+})}
\le
\left\| \partial_{\bf I}f-\mathfrak{S}_hf\right\|_{\mathcal{B}(\RR^d_{m,+})}+\| \mathfrak{S}_h\|\|f\|_{\mathcal{B}(\RR^d_{m,+})}
\\ \le 
\frac {\| \partial_
{\bf{I}}f\|_{H^\omega(\RR^d_{m,+})} }{2^{d-m}h^d}\int_{B_h}\omega(\rho (u,\theta))du+\frac 1{2^{d-m}h^d}\sup_{x\in \RR^d_{m,+}}\left|\int_{x+B_h}\partial_
{\bf{I}}f(u)du\right|.
\end{multline*}
Representation~\eqref{mixdiff} implies
\[
\frac 1{2^{d-m}h^d}\sup_{x\in \RR^d_{m,+}}\left|\int_{x+B_h}\partial_
{\bf{I}}f(u)du\right|\le \frac {2^d}{2^{d-m}h^d}\| f\|_{\mathcal{B}(\RR^d_{m,+})}=\frac {2^m}{h^d}\| f\|_{\mathcal{B}(\RR^d_{m,+})},
\]
which implies inequality~\eqref{partial_omega}.

In the case $\omega(t) = t^\alpha$, using  the layer cake representation, see e.g.,~\cite[Theorem 1.13]{lieb2001}, symmetry considerations, and writing $|x|_\infty$ instead of $\rho(x,\theta)$, we obtain 
\begin{multline*}
\int_{B_h}\omega(\rho (u,\theta))du 
=
2^{d-m}\int_{(0,h)^d}|u|_\infty^\alpha du
\\=
2^{d-m}\int_0^\infty \mu\left\{v\in (0,h)^d\colon |u|_\infty^\alpha > t\right\}dt 
=
2^{d-m}\int_0^{h^{\alpha}} \left(h^d - t^{\frac d\alpha}\right)dt 
\\= 
2^{d-m}h^{d+\alpha}\left(1 - \frac{\alpha}{d+\alpha}\right)
= \frac{d\cdot 2^{d-m}}{d+\alpha}h^{d+\alpha}. 
\end{multline*}
So that the right-hand side of~\eqref{partial_omega} becomes
$$
\frac {d}{d+\alpha}\| \partial_
{\bf{I}}f\|_{H^\omega(\RR^d_{m,+})}h^\alpha +\frac {2^m}{h^d}\| f\|_{\mathcal{B}(\RR^d_{m,+})}.
$$
Minimizing this expression with respect to $h>0$ i.e., choosing 
$$
h = 2^{\frac{m}{d+\alpha}}\left(\frac{d+\alpha}{\alpha}\cdot\frac{\| f\|_{\mathcal{B}(\RR^d_{m,+})}}{\| \partial_
{\bf{I}}f\|_{H^\omega(\RR^d_{m,+})}}\right)^{\frac 1{d+\alpha}}
$$
we obtain the right-hand side of~\eqref{partialOmegaMultiplicative}.

Next we prove sharpness of inequality~\eqref{partial_omega} for the case $m=0$. Consider the function
$$
g_{e,h}(x) = \int_0^{x_1}\ldots\int_0^{x_d}f_{e,h}(u)du, 
$$
where $f_{e,h}$ is defined in~\eqref{f_eh}. Then  $\partial_{\bf I}g_{e,h} = f_{e,h}$, and hence 
\[
\| \partial_{\bf I}g_{e,h}\|_{\mathcal{B}(\RR^d)}=\omega(h)\text{ and }
\| \partial_{\bf I}g_{e,h}\|_{H^\omega(\RR^d)}=1.
\]
Moreover,
\[
\| g_{e,h}\|_{\mathcal{B}(\RR^d)}
=
\int_0^{h}\ldots\int_0^{h}(\omega(h)-\omega(|u|_\infty))du
=
h^d\omega(h) -\int_0^{h}\ldots\int_0^{h}\omega(|u|_\infty)du
\]
and due to symmetricity of $[-h,h]^d$,
\[
\int_{B_h}\omega(|u|_\infty)du
=
2^d\int_0^{h}\ldots\int_0^{h}\omega\left(|u|_\infty\right)du.
\]
Direct computations now show that inequalities~\eqref{partial_omega} and~\eqref{partialOmegaMultiplicative} become equality for the function $g_{e,h}$.

Finally, we prove sharpness of inequality~\eqref{partial_omega} in the case $m = 1$. In this case $B_h=(0,h)\times (-h,h)^{d-1}$.
There exists $0<a<h$ such that 
\begin{multline*}
\int_{\{ x\in B_h\colon x_1<a\}}(\omega(h)-\omega(|x|_\infty))dx=\int_{\{ x\in B_h\colon x_1>a\}}(\omega(h)-\omega(|x|_\infty))dx
\\=
\frac 12\int_{B_h}(\omega(h)-\omega(|x|_\infty))dx.
\end{multline*}
The set $B_h$ consists of $2^{d-1}$ equal cubes with edge lengths equal to $h$; $\theta =\theta_d$ is one of the vertices for each of these cubes. The hyperplane  $x_1=a$ divides these cubes into pieces  $c_1^-,\ldots ,c_{2^{d-1}}^-$ that have  $\theta$ among vertices, and pieces $c_1^+,\ldots ,c_{2^{d-1}}^+$ that have $(a,\theta_{d-1})$ among their vertices. 

It is clear that 
\[
\int_{\{ x\in B_h\colon x_1<a\}}(\omega(h)-\omega(|x|_\infty))dx=\sum_{i=1}^{2^{d-1}}\int_{c_i^-}(\omega(h)-\omega(|x|_\infty))dx
\]
and
\[
\int_{\{ x\in B_h\colon x_1>a\}}(\omega(h)-\omega(|x|_\infty))dx=\sum_{i=1}^{2^{d-1}}\int_{c_i^+}(\omega(h)-\omega(|x|_\infty))dx.
\]
From the symmetry considerations it follows that each of $2^{d-1}$ summands in each of the right-hand sides of these equalities are equal. Thus for each $i,j=1,\ldots,2^{d-1}$
\[
\int\limits_{c_i^-}(\omega(h)-\omega(|x|_\infty))dx=\int\limits_{c_j^+}(\omega(h)-\omega(|x|_\infty))dx=\frac 1{2^d}\int\limits_{B_h}(\omega(h)-\omega(|x|_\infty))dx.
\]
An extremal function in this case we define as follows.
\[ G_{e,h}(x)=\int_a^{x_1}\int_0^{x_2}\ldots\int_0^{x_d}\partial_{\bf I}f_{e,h}(u)du_d\ldots du_1,
 \]
 where $f_{e,h}$ is defined in~\eqref{f_eh}. For this function we have
 $$
 \| G_{e,h}\|_{\mathcal{B}(\RR^d_{1,+})}=\frac 1{2^d}\int_{B_h}(\omega(h)-\omega(|x|_\infty))dx
 =
\frac{h^d}{2}\omega(h) - \frac{1}{2^d}\int_{B_h}\omega(|x|_\infty)dx, 
 $$
 $$
 \| \partial_{\bf I}G_{e,h}\|_{\mathcal{B}(\RR^d_{1,+})}=\omega(h),\text{ and }\| \partial_{\bf I}G_{e,h}\|_{H^\omega(\RR^d_{1,+})}=1.
$$
Direct computations now show that inequalities~\eqref{partial_omega} and~\eqref{partialOmegaMultiplicative} with $m=1$ become equalities on the function $G_{e,h}$.
\end{proof}

\begin{remark}
It is not clear, whether inequality~\eqref{partial_omega} is sharp for $m=2,\ldots,d$, even in the case $d=2$ and $\omega(t)=t$.
\end{remark}

\bibliographystyle{elsarticle-num}
\bibliography{mybibfile}

\begin{thebibliography}{10}
\expandafter\ifx\csname url\endcsname\relax
  \def\url#1{\texttt{#1}}\fi
\expandafter\ifx\csname urlprefix\endcsname\relax\def\urlprefix{URL }\fi
\expandafter\ifx\csname href\endcsname\relax
  \def\href#1#2{#2} \def\path#1{#1}\fi

\bibitem{Landau}
E.~Landau, {\"U}ber einen satz des herrn {L}ittlewood, Rend. Circ. Mat. Palermo
  (1884-1940) 35 (1913) 265--276.

\bibitem{Kolmogorov1939}
A.~Kolmogorov, On inequalities between the upper bounds of the successive
  derivatives of an arbitrary function on an infinite interval, Uchenye Zapiski
  MGU. Math 30~(3) (1939) 3--13, (in Russian).

\bibitem{Nagy}
B.~Sz.-Nagy, {\"U}ber {I}ntegralungleichungen zwischen einer {F}unction und
  ihrer {A}bleitung, Acta Sci. Math. 10 (1941) 64--74.

\bibitem{Arestov1996}
V.~Arestov, Approximation of unbounded operators by bounded operators and
  related extremal problems, Russ. Math. Surv. 51~(6) (1996) 1093.
\newblock \href {https://doi.org/10.1070/RM1996v051n06ABEH003001}
  {\path{doi:10.1070/RM1996v051n06ABEH003001}}.

\bibitem{BKKP}
V.~Babenko, N.~Korneichuk, V.~Kofanov, S.~Pichugov, Inequalities for
  derivatives and their applications., Naukova Dumka, Kiev, 2003, (in Russian).

\bibitem{Babenko2012}
V.~F. Babenko, Inequalities of Kolmogorov type for fractional derivatives and
  their applications, Problems and methods: Mathematics, Mechanics,
  Cibernetics, Vol. 4, Chapter 2, Naukova knyga, 2012, (in Russian).

\bibitem{Babenko22}
V.~Babenko, O.~Kovalenko, N.~Parfinovych, On approximation of hypersingular
  integral operators by bounded ones, J. Math. Anal. Appl. 513~(2) (2022)
  126215.
\newblock \href {https://doi.org/10.1016/j.jmaa.2022.126215}
  {\path{doi:10.1016/j.jmaa.2022.126215}}.

\bibitem{Babenko05}
V.~Babenko, V.~Kofanov, S.~Pichugov, Inequalities of {N}agy type for periodic
  functions, East J. Approx. 11~(1) (2005) 3--11.

\bibitem{Kofanov}
V.~Kofanov, I.~Popovich, Sharp {N}agy type inequalities for the classes of
  functions with given quotient of the uniform norms of positive and negative
  parts of a function, Res. Math. 28~(1) (2020) 3--11.
\newblock \href {https://doi.org/10.15421/242001} {\path{doi:10.15421/242001}}.

\bibitem{Babenko23}
V.~F. Babenko, V.~V. Babenko, O.~V. Kovalenko, N.~V. Parfinovych, On {L}andau
  – {K}olmogorov type inequalities for charges and their applications, Res.
  Math. 31~(1) (2023) 3--16.
\newblock \href {https://doi.org/10.15421/242301} {\path{doi:10.15421/242301}}.

\bibitem{Nikolsky46}
S.~M. Nikol'skii, Fourier series of functions with a given modulus of
  continuity, Dokl. Akad. Nauk SSSR 52~(3) (1946) 191--194.

\bibitem{ExactConstants}
N.~Korneichuk, Exact Constants in Approximation Theory, Encyclopedia of
  Mathematics and its Applications, Cambridge University Press, 1991.

\bibitem{Babenko16}
V.~Babenko, Y.~Babenko, N.~Parfinovych, D.~Skorokhodov, Optimal recovery of
  integral operators and its applications, J. Complexity 35 (2016) 102--123.

\bibitem{VeraBabenko_JANO}
V.~F. Babenko, V.~V. Babenko, Best approximation, optimal recovery, and
  {L}andau inequalities for derivatives of {H}ukuhara-type in function
  {L}-spaces, J. Appl. Numer. Optim. 1 (2019) 167--182.

\bibitem{Kovalenko20}
O.~Kovalenko, On optimal recovery of integrals of random processes, J. Math.
  Anal. Appl. 487~(1) (2020) 123949.

\bibitem{Babenko21}
V.~Babenko, V.~Babenko, O.~Kovalenko, M.~Polishchuk, Optimal recovery of
  operators in function {L}-spaces., Anal. Math. 47 (2021) 13–32.
\newblock \href {https://doi.org/10.1007/s10476-021-0065-y}
  {\path{doi:10.1007/s10476-021-0065-y}}.

\bibitem{Kriachko}
N.~Kriachko, Two sharp inequalities for operators in a {H}ilbert space, Res.
  Math. 30~(1) (2022) 56--65.
\newblock \href {https://doi.org/10.15421/242206} {\path{doi:10.15421/242206}}.

\bibitem{BabenkoArxiv}
V.~Babenko, V.~V. Babenko, O.~Kovalenko, Korneichuk-{S}techkin lemma,
  {O}strowski and {L}andau inequalities, and optimal recovery problems for
  ${L}$-space valued functions, arXiv:2006.14581 (2021) 1--25\href
  {https://doi.org/10.48550/arXiv.2006.14581}
  {\path{doi:10.48550/arXiv.2006.14581}}.

\bibitem{Babenko22b}
V.~Babenko, O.~Kovalenko, N.~Parfinovych, General form of
  $(\lambda,\varphi)$--additive operators on spaces of ${L}$-space-valued
  functions, Res. Math. 30~(1) (2022) 3--9.

\bibitem{Babenko07}
V.~F. Babenko, M.~S. Churilova, Kolmogorov type inequalities for hypersingular
  integrals with homogeneous characteristic, Banach J. Math. Anal. 1~(1) (2007)
  66 -- 77.
\newblock \href {https://doi.org/10.15352/bjma/1240321556}
  {\path{doi:10.15352/bjma/1240321556}}.

\bibitem{babenko2010}
V.~F. Babenko, D.~A. Levchenko, Kolmogorov type inequalities for hypersingular
  integrals with sign-alternating characteristic, Res. Math. 15 (2010) 18--27,
  (in Russian).

\bibitem{BPP2010}
V.~Babenko, N.~Parfinovich, S.~Pichugov, Sharp {K}olmogorov-type inequalities
  for norms of fractional derivatives of multivariate functions, Ukr. Math. J.
  62 (2010) 343--357.
\newblock \href {https://doi.org/10.1007/s11253-010-0358-y}
  {\path{doi:10.1007/s11253-010-0358-y}}.

\bibitem{Parfinovych}
N.~Parfinovych, V.~Pylypenko, Kolmogorov inequalities for norms of
  marchaud-type fractional derivatives of multivariate functions, Res. Math.
  28~(1) (2020) 10--23.
\newblock \href {https://doi.org/0.15421/242007} {\path{doi:0.15421/242007}}.

\bibitem{Stechkin1967}
S.~Stechkin, Best approximation of linear operators, Math. Notes 1~(2) (1967)
  91--99.

\bibitem{Loomis}
L.~Loomis, An introduction to abstract harmonic analysis, D. Van Nostrand
  Company, 1953.

\bibitem{ambrosio2004topics}
L.~Ambrosio, P.~Tilli, Topics on Analysis in Metric Spaces, Volume 25 of Oxford
  Lecture Mathematics, Oxford University Press, 2004.

\bibitem{heinonen2015sobolev}
J.~Heinonen, P.~Koskela, N.~Shanmugalingam, J.~Tyson, Sobolev Spaces on Metric
  Measure Spaces: An Approach Based on Upper Gradients, New Mathematical
  Monographs, Cambridge University Press, 2015.

\bibitem{Berezanski}
Y.~Berezanski, G.~Us, Z.~G. Sheftel, Functional analysis, Elsevier Science,
  2003.

\bibitem{lieb2001}
E.~Lieb, M.~Loss, Analysis, Crm Proceedings \& Lecture Notes, American
  Mathematical Society, 2001.

\end{thebibliography}

\end{document}